\documentclass[conference]{IEEEtran}
\IEEEoverridecommandlockouts
\usepackage{graphicx}
\usepackage{hyperref}
\usepackage{pstricks}
\usepackage{pst-node}
\usepackage{pst-tree}
\usepackage[font=small]{caption}
\usepackage{subcaption}
\usepackage{lingmacros}
\usepackage{tikz}
\usepackage{pgf}
\usepackage{pgffor}
\usepackage{pgfmath}
\usepackage{amsmath}
\usepackage{cite}
\usepackage{amsmath,amssymb,amsfonts}
\usepackage{algorithmic}
\usepackage{graphicx}
\usepackage{textcomp}
\usepackage{xcolor}
\usepackage{array}
\usepackage{makecell}
\usepackage{cite}
\def\BibTeX{{\rm B\kern-.05em{\sc i\kern-.025em b}\kern-.08em
    T\kern-.1667em\lower.7ex\hbox{E}\kern-.125emX}}
\begin{document}

\title{Dynamic Network Flow Optimization for Task Scheduling in PTZ Camera Surveillance Systems\\
\thanks{This material is based upon work supported by the U.S. Department of Homeland Security under [Grant Award 22STESE00001-02-00] made through Northeastern University. The views and conclusions contained in his document are those of the authors and should not be interpreted as necessarily representing the official policies, either expressed or implied, of the U.S. Department of Homeland Security.}
}

\author{\IEEEauthorblockN{1\textsuperscript{st} Mohammad Merati*}
\IEEEauthorblockA{\textit{Division of Systems Engineering} \\
\textit{Boston University}\\
Boston, USA \\
mmerati@bu.edu}
*Corresponding author
~\\
\and
\IEEEauthorblockN{2\textsuperscript{nd} David Castañón}
\IEEEauthorblockA{\textit{Department of Electrical and Computer Engineering} \\
\textit{Boston University}\\
Boston, USA \\
dac@bu.edu}
}

\maketitle

% ----------  mandatory IEEE notice  ----------
\begingroup
\renewcommand\thefootnote{}\phantomsection\footnotetext{%
\textcopyright~2025 IEEE. Personal use of this material is permitted. Permission from IEEE must be obtained for all other uses, in any current or future media, including reprinting/republishing this material for advertising or promotional purposes, creating new collective works, for resale or redistribution to servers or lists, or reuse of any copyrighted component of this work in other works.\\}
\endgroup
% ---------------------------------------------

\begin{abstract}
This paper presents a novel approach for optimizing the scheduling and control of Pan-Tilt-Zoom (PTZ) cameras in dynamic surveillance environments. The proposed method integrates Kalman filters for motion prediction with a dynamic network flow model to enhance real-time video capture efficiency. By assigning Kalman filters to tracked objects, the system predicts future locations, enabling precise scheduling of camera tasks. This prediction-driven approach is formulated as a network flow optimization, ensuring scalability and adaptability to various surveillance scenarios. To further reduce redundant monitoring, we also incorporate group-tracking nodes, allowing multiple objects to be captured within a single camera focus when appropriate.
In addition, a value-based system is introduced to prioritize camera actions, focusing on the timely capture of critical events. By adjusting the decay rates of these values over time, the system ensures prompt responses to tasks with imminent deadlines. Extensive simulations demonstrate that this approach improves coverage, reduces average wait times, and minimizes missed events compared to traditional master-slave camera systems. Overall, our method significantly enhances the efficiency, scalability, and effectiveness of surveillance systems, particularly in dynamic and crowded environments.
\end{abstract}


\section{Introduction}
The landscape of surveillance technology has rapidly evolved, leading to the widespread adoption of smart video capturing systems. Among these, Pan-Tilt-Zoom (PTZ) cameras have become essential tools due to their ability to dynamically adjust their viewing angles and zoom levels. This flexibility allows for effective monitoring of large areas with fewer cameras, making PTZ cameras a critical component in applications such as security surveillance, traffic monitoring, and public safety.
PTZ cameras' primary advantage lies in their ability to focus on specific targets by panning, tilting, and zooming. This adaptability is beneficial in environments where the movement and behavior of objects need to be closely monitored in real time. Managing these cameras efficiently presents significant challenges. Effective control algorithms are needed to predict object movements, schedule camera tasks, and optimize coverage.

A key technique employed in this domain is the use of Kalman filters \cite{welch1995introduction} for motion prediction. These filters help in tracking objects and predicting their future locations, which is crucial for determining where and when a PTZ camera should focus. To complement this, one needs scheduling algorithms that can make quick and efficient decisions to manage the camera's tasks in real-time. 

There has been much past work on high-speed object tracking and image quality enhancement using control real-time pan-tilt adjustments, zoomed-in tracking, and depth-of-field extension.  Some approaches focus on using specialized hardware without addressing broader task scheduling and surveillance optimization needed in dynamic environments \cite{fast1,fast2,fast3}. Other approaches use the concept of Task Visibility Intervals (TVIs) to optimize PTZ camera scheduling.  The authors use a greedy algorithm that focusing on motion prediction and occlusion  without employing a Kalman filter for tracking \cite{TVI1,TVI2}.

The most common approach in the literature is to use a master-slave system where fixed cameras are used for wide-area detection and tracking of people, while PTZ cameras are controlled to perform close-in inspection on individual persons. \cite{Bagdanov} focuses on high-resolution image capture using a master-slave system. Reference \cite{Xu} presents an autonomous vision system using multiple robotic PTZ cameras and a fixed wide-angle camera. It employs a greedy scheduling approach to maximize request satisfaction. Reference \cite{SurveillanceCameraScheduling} formulates camera control as an online scheduling problem and uses a weighted round-robin heuristic to guide PTZ cameras for pedestrian tracking.
\cite{proenca} discusses a system for biometric recognition using an architecture with a wide-view static camera and a PTZ camera which is used to capture high-resolution data on the head regions of individuals. 

Reference \cite{cai} uses a single PTZ camera to track individuals and capture high-resolution facial images with a scheduling algorithm that zooms in based on detection order.
For multi-camera coordination \cite{Segawa} discusses strategies for managing a network of PTZ cameras that work together to track multiple targets, without having to schedule zoomed-in views.
Five scheduling policies are compared in \cite{Ward} to optimize real-time camera assignments, focusing on maximizing coverage without using zooming for additional information.
Similarly, \cite{Neves} addresses the time challenges of dynamic camera scheduling for tracking persons in crowded scenes.

Reference \cite{Ding} optimizes distributed camera networks to capture high-resolution images of individual and group activities. It employs local decision-making for camera coordination and emphasizes opportunistic image acquisition during events, enhancing image quality through distributed optimization. A reactive coverage control algorithm for PTZ cameras is presented in \cite{Arslan}, using conic Voronoi diagrams to optimize orientations and zoom levels while adjusting focus based on spatial event distributions.
Reference \cite{Bansal2011} proposed a heuristic for load balancing in distributed systems, and  \cite{Gawali2017} introduced an optimization algorithm for efficient task scheduling in cloud computing. However, both approaches are limited to static or semi-dynamic environments and lack real-time adaptability.
A system for controlling PTZ cameras to capture high-quality imagery for biometric tasks is presented in \cite{Krahnstoever}, which collaborates with fixed cameras for multi-target tracking and uses graph-based greedy heuristics for efficient scheduling.

Despite these advances, existing methods often struggle with scalability and efficiency in highly dynamic settings, and cannot handle the computation complexity of scheduling PTZ cameras across both tracking and high-resolution imaging.  There is a need for a robust solution that can adapt to varying numbers of targets and complex movement patterns while ensuring that critical events are captured promptly and with high quality.

In this paper, we present a novel approach to enhance the scheduling and control of PTZ cameras in video surveillance networks. Our method integrates Kalman filters for motion prediction with a set cover algorithm to optimize camera scheduling. By assigning a Kalman filter to each tracked object, we accurately predict future locations, which is crucial for generating tracking nodes. These nodes are used to construct a dynamic network flow model, allowing us to formulate the scheduling problem as a network flow optimization that can be solved efficiently in real-time.

A key innovation of our approach is the introduction of group-tracking nodes. Instead of treating each target individually, we group targets based on their predicted locations, significantly reducing the number of groups the system needs to monitor. This minimizes redundant tracking, where multiple cameras might otherwise capture the same target repeatedly, and allows for more efficient utilization of camera resources.

This paper addresses key challenges in PTZ camera scheduling and control for joint tracking and high-resolution video capturing by proposing a robust solution using advanced prediction and optimization techniques. It presents concepts to enhance the efficiency and reliability of smart video surveillance, tailored to the demands of complex, dynamic environments. Through simulations, we demonstrate that our method outperforms existing algorithms, capturing more critical events and providing higher-quality footage.
\section{Problem Formulation}
In this section, we present a dynamic network flow model formulated as an integer optimization problem for online video surveillance. The model includes several key components: the dynamic network itself, flow variables, nodes, arcs, and time-dependent constraints.

\subsection{Main Approach}

We begin by describing the dynamic network that represents the surveillance system. This network is modeled as a graph where nodes represent both cameras and locations within the scene. Each node is specified by a tuple $(i, t) \in N(t)$, where $i$ denotes the node identifier and $t$ represents the time step. The set $N(t)$ contains all nodes existing at time step $t$. There is a node $R_{i,t}$ for each camera at time $t$. There are location nodes of two types
\begin{enumerate}
    \item Track Location Node $K_{j,t}$ represents the predicted location of individual object $j$ at time $t$, where a camera may zoom in to focus on that object. 
    \item Fixed Location Node $F_{j,t}$ represents a region that can be observed in wide area mode at time $t$, collecting low resolution video on all objects in that region.  These regions are predefined in the scene, and depend on camera locations. 
    \end{enumerate}
  
In our planning problem, we want to visit each fixed location node at least once every $T$ time periods to guarantee that we maintain tracks on all objects in a region.  We introduce accumulation nodes $D_(j,\tau)$ for every fixed location $j$ at times $t = \tau T$ that accumulate how many times location $j$ has been visited by a camera node over times $t = (\tau-1)T + 1, \ldots, \tau T$.  

Arcs in the network connect camera nodes to location nodes or location nodes to other location nodes. An arc from a camera node $R_i$ to location node $K_j$ or $F_j$ at time $t$ indicates this camera can observe that location at time $t$.  An arc from location node $K_j$ at time $t$ to location node $K_j$ at time $t'$ is used to remember whether a location has already been observed in high resolution.  An arc from location node $F_j$ at time $t$ to location node $D_j$ at time $t = \tau T$ is used to whether this fixed location has been observed already in the past $T$ periods.  

For each time step $t$, let $A(t)$ denote the set of arcs originating from nodes in $N(t)$. Each arc $((i, t), (j, t')) \in A(t)$ can carry a flow that represents the assignment of a camera $i$ to a location $j$ at a time $t$.  The entire set of arcs in the problem is given by $\mathcal{E} = A(1) \cup A(2) \cup \dots \cup A(T)$.  The entire set of nodes in the problem is defined as $\mathcal{V} = N(1) \cup N(2) \cup \dots \cup N(T)$, where $T$ is the total number of time steps.

The network $\mathcal{G} = (\mathcal{V}, \mathcal{E})$ models the possible assignments of cameras to locations over time. In this network:
\begin{itemize}
    \item The camera nodes act as supply nodes, each with a supply amount of one at each time, meaning each camera can be assigned to only one action at each time step.
    \item A sink node $S$ is created to accumulate the total camera views over time from both individual objects as well as regions observed.
    \item Demand nodes $D_{j,\tau}$ have a demand of 1 to require at least one look at each fixed location every $T$ periods.
    \item The capacities of the arcs are set to one, except for those connecting demand $D_{j,\tau}$ to sink nodes to send excess looks over a window of $T$ periods to the sink $S$.
\end{itemize}

\begin{figure}
    \centering
    \includegraphics[scale=0.45]{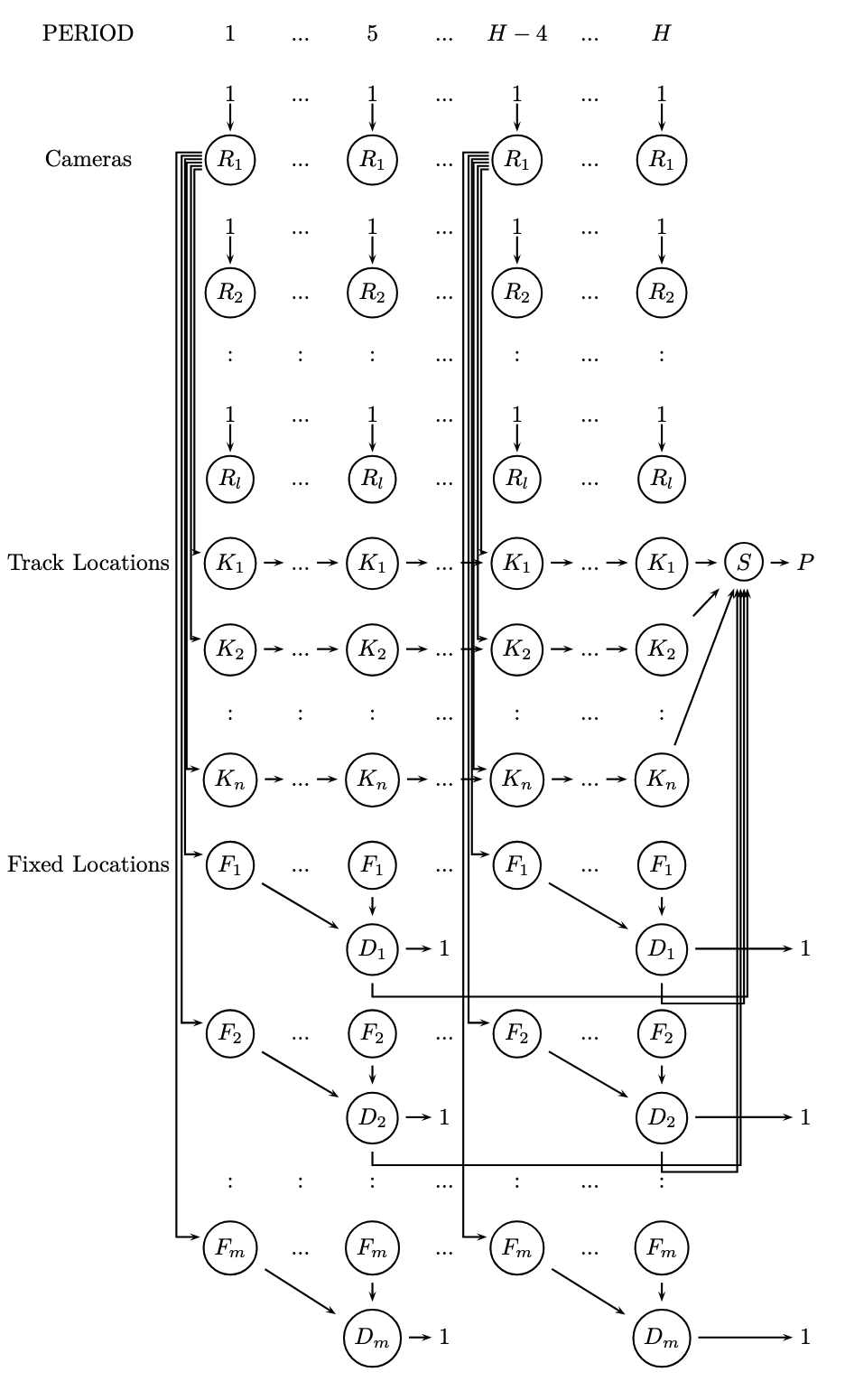}
    \caption{ H-period online video surveillance network model with $T=5$. The graph includes $l$ camera nodes, $n$ track location nodes, and $m$ fixed location nodes, each duplicated across H periods. Node $D_i$ denotes the demand node for fixed locations which are repeated $H/T$ times. For simplicity, the $t$ is omitted in designating node $(X_i,t)$. }
    \label{fig:graph}
    \vspace{-10pt}
\end{figure}

Fig. \ref{fig:graph} illustrates our model for a system with $l$ cameras, $n$ track locations (moving objects), and $m$ fixed locations. Each row in the figure represents a single entity (camera or location) repeated over time. 
With the network components defined, we can formulate the dynamic network flow problem as follows:

\begin{small} 
\begin{align} \label{objAround}
    \text{maximize} \quad &\sum_{t=1}^{H}\sum_{i=1}^{l} \bigg(\sum_{j=1}^{n} c_{(R_i,t)(K_j,t)} x_{(R_i,t)(K_j,t)} \\ 
    & + \sum_{k=1}^{m} c_{(R_i,t)(F_k,t)} x_{(R_i,t)(F_k,t)} \bigg)\nonumber 
\end{align}
\vspace{-15pt}
\begin{align} \label{const1}
    \text{s.t.} \quad \sum_{j=1}^{n} {x_{(R_i,t)(K_j,t)}} = 1 \quad i=1,...,l 
    ,\ t=1,...,H your
\end{align}
\vspace{-15pt}
\begin{align} \label{const2}
    \sum_{i=1}^{l} &{x_{(R_i,t)(K_j,t)}} + x_{(K_j,t-1)(K_j,t)} = x_{(K_j,t)(K_j,t+1)}\ \  \forall\ j, t
\end{align}
\vspace{-15pt}
\begin{align} \label{const3}
    \sum_{i=1}^{l} {x_{(R_i,H)(K_j,H)}} + x_{(K_j,H-1)(K_j,H)} = x_{(K_j,H)S}\quad \forall\ j
\end{align}
\vspace{-15pt}
\begin{align} \label{const4} 
    \sum_{i=1}^{l} {x_{(R_i,t)(F_j,t)}} = x_{(F_j,t)(D_j,\tau)} \quad  \forall\ j, \ t, \  \tau =\lceil t/T \rceil
\end{align}
\vspace{-15pt}
\begin{align} \label{const5}
    \sum_{k=0}^{T-1} {x_{(F_j,T\tau-k)(D_j,\tau)}} = 1 + x_{(D_j,\tau)S} \quad  \forall\ j, \tau = 1, \ldots, \frac{H}{T}
\end{align}
\vspace{-15pt}
\begin{align} \label{const6}
    \sum_{i=1}^{n} {x_{(K_i,H)S}} + \sum_{\tau=1}^{H/T}{\sum_{j=1}^{m} {x_{(D_j,\tau)S}}} = P
\end{align}
\vspace{-15pt}
%\begin{flalign} \label{const7}
%    p_{(R_i,t)(K_j,t)} x_{(R_i,t)(K_j,t)} \leq e_{(K_j,t)} \quad & \forall\ i,\ j,\ t 
%%\end{flalign}
%\vspace{-10pt}
%\begin{flalign} \label{const8}
%    x_{(R_i,t)(K_j,t)} \leq y_{(R_i,t)(K_j,t)} \quad \forall\ i,\ j,\ t 
%\end{flalign}
%\vspace{-10pt}
\begin{flalign} \label{const9}
    x_{(R_i,t)(K_j,t)} \in \{0,1\} \quad & \forall\ i,\ j,\ t 
\end{flalign}
\vspace{-15pt}
\begin{flalign} \label{const10}
    x_{(K_j,t)(K_j,t+1)} \in \{0,1\} \quad & \forall\  j,\ t 
\end{flalign}
\vspace{-15pt}
\begin{flalign} \label{const11}
    x_{(R_i,t)(F_j,t)} \in \{0,1\} \quad & \forall\ i,\ j,\ t 
\end{flalign}
\vspace{-15pt}
\begin{flalign} \label{const12}
    x_{(F_j,T\tau-k)(D_j,\tau)} \in \{0,1\} \quad & \forall \ j,\ k \in 0, \ldots,T-1  
\end{flalign}
\vspace{-15pt}
\begin{flalign} \label{const14}
    x_{(K_i,H)S} \in \{0,1\} \quad i=1,...,n
\end{flalign}
\vspace{-15pt}
\begin{flalign} \label{const15}
    x_{(D_j,\tau)S} \in \{0, \ldots, T-1\} \quad \forall\ j, \  \tau=1,...,H/T 
\end{flalign}
\end{small} 
In this formulation, we aim to maximize objective function \eqref{obj} which is the total value of camera assignments to locations over the planning horizon $H$. The coefficients $c_{(R_i,t)(K_j,t)}$ and $c_{(R_i,t)(F_k,t)}$ represent the values of assigning camera $R_i$ to track location $K_j$ and fixed location $F_k$ at time $t$, respectively. Variables $x_{(R_i,t)(K_j,t)}$ and $x_{(R_i,t)(F_k,t)}$ are flows originating at camera node $R_i$ to track node $K_j$ and fixed location $F_k$ at period $t$, respectively. Equations \eqref{const1} ensures that each camera is assigned to exactly one action (either a track location or a fixed location) at each time step. Equations \eqref{const2} and \eqref{const3} represent flow conservation for track location nodes, ensuring that each track is viewed at high resolution at most once in the planning period.   Equations \eqref{const4} and \eqref{const5} handle flow conservation for fixed location nodes, ensuring that fixed locations are monitored at least once every $T$ time steps. Variable $x_{(F_j,t)(D_j,\tau)}$ is the flow from fixed location node $F_j$ at time $t$ to the demand node $D_j$, and $x_{(D_j,\tau)S}$ is the flow from demand node $D_{j,\tau}$ to sink node $S$. Equation \eqref{const6} is the constraint to balance the total number of assignments with the available opportunities for the cameras, accounting for the periodic monitoring of fixed locations. The parameter $P$ represents  total  opportunities and is calculated as $
    P = T \times l - (H/T)\times m. $
The other constraints indicate that decision variables are 0-1.

\subsection{Group-tracking Task Assignment}
To enhance the efficiency of our approach, we introduce a modification that leverages group tracking. The rationale behind this modification is based on the observation that pedestrians often move in close proximity within crowded environments. This proximity allows for the possibility of capturing multiple objects with a single camera high resolution interrogation, thereby optimizing camera usage and reducing redundant tracking.

In our model, we utilize Kalman filters to predict the future locations of each tracked object. By assigning a Kalman filter to every object, we can accurately estimate its future position based on its current state and motion dynamics. At each decision period, we generate predicted tracking nodes for all active tracks, representing the expected positions of the objects in upcoming time steps.

Using one camera as a reference, we calculate the coverage areas that would result from focusing on the center of each predicted tracking node. For each coverage area, we identify the set of track IDs (objects) that would be within view if the camera focuses on that area.  To efficiently select a minimal set of coverage areas that collectively encompass all active tracks, we apply a greedy set cover algorithm. This algorithm prioritizes coverage areas that include the largest number of uncovered tracks, iteratively adding them to the solution until all tracks are covered. By doing so, we form group tracking nodes where each node represents a cluster of tracks that can be monitored simultaneously by focusing on a specific area.

By grouping the targets in this manner, we effectively reduce the number of tracking nodes the system needs to consider. This not only improves overall performance by decreasing computational complexity but also minimizes the likelihood of tracking the same object multiple times. Consequently, the system has more opportunities to capture new objects, enhancing the efficiency of surveillance.

The introduction of group tracking modifies the network flow model by replacing individual track location nodes with group tracking nodes. The flow conservation and assignment constraints are adjusted accordingly to account for the grouped targets, ensuring that the optimization problem remains tractable and can be solved efficiently. This modification aligns with our goal of optimizing camera scheduling while maintaining high-quality coverage of all targets in the scene.

\section{Value System}
In this section, we explain how to assign values to the arcs connecting cameras to locations in our network model. First we assign values to the arcs connecting cameras to fixed-location nodes. We use two main characteristics:
\begin{itemize}
\item the number of tracks $N_{F_j,t}$ that would be within the field of view  if camera $i$ focuses on the fixed location. 
\item The angle between the camera's line of sight towards the fixed location and a reference direction, denoted as $e_{R_i,F_j}(t)$. 
\end{itemize}
We compute the video quality value $V(e_{R_i,F_j}(t))$ as
$$ V(e_{R_i,K_j}(t)) = \begin{cases} 3 & |e_{R_i,F_j}(t)| \le \pi/6\\ 2 & |e_{R_i,K_j}(t)| \in (\pi/6, \pi/3]\\ 1 & |e_{R_i,K_j}(t)| \in (\pi/3, \pi/2]\\ 0 & \text{elsewhere} \end{cases} $$ 

The formula is $c_{(R_i,t)(F_j,t)} = N_{F_j,t} + V(e_{R_i,K_j}(t))$.

Second, we assign values to the arcs connecting cameras to individual cluster locations.  
We use the following information: 
\begin{itemize}
\item The predicted exit time for cluster $K_j$ from the scene $t^e_{K_j}$.
\item The size $N_{K_j}$ of the cluster of tracks that can be observed simultaneously.
\item The alignment between the camera's observation direction and the track's movement direction as an angle $e_{R_i,K_j}(t)$. 
\end{itemize}

We define the number of tracks departing in planning period as $N^e = \sum_{j =1}^n I\{t^e_j < H\}$.  Then, the number 
that are not departing during the planning period is $N^s = n - N^e$.  For tracks that are not departing, we assign values as 
\begin{align*}
c_{(R_i,t)(F_j,t)} &= [(N^s + 1)(H - (t-1)) + 
 N^s \\&- rank(t^e_{K_j})+ V(e_{R_i,K_j}(t))]\times N_{K_j}
\end{align*} 
where $rank(t^e_{K_j})$ is the rank in terms of increasing departure time for the tracks that leave the area of surveillance after the planning horizon.

For tracks that are departing in the planning period, we assign to assign higher values.  Let $V^{1} = (N^s + 1)H$. Then, 
$$  c_{(R_i,t)(F_j,t)} = [V^1 \times 2^{N^e + 1 - rank(t^e(K_j))} +  V(e_{R_i,K_j}(t))]\times N_{K_j}$$
where $rank(t^e(K_j))$ is the order in earliest time of departure for tracks $K_j$ that are departing in this planning period.  In this manner, tracks which are departing the surveillance area in the near future get higher values, and are thus more likely to be inspected at high resolution before leaving the area.

\section{Experimental Study}
In this section, we present an empirical evaluation of our proposed network flow model and compare its performance with that of a master-slave camera system. Both systems utilize the same number of cameras; however, in the master-slave configuration, one of the PTZ cameras is replaced with a wide-view static camera. We begin by describing the simulation scenarios in section \ref{sec:simulation-scenario}, we describe the simulated scenarios, outline the metrics used for evaluation in section \ref{sec:metrics} and then discuss the results in section \ref{sec:resutls}.
\subsection{Simulation Scenarios} \label{sec:simulation-scenario}
To assess the effectiveness of our surveillance system, we designed simulation scenarios that mimic real-world crowded environments. Pedestrians are introduced into the environment based on a Poisson process, entering from the upper side of the simulated area. Their movements follow a constant velocity model within a two-dimensional Cartesian coordinate system.

In the flexible system configuration, the environment is monitored using three PTZ cameras positioned along the bottom side. This placement ensures that pedestrians entering the scene face towards the cameras, eliminating issues related to capturing faces due to unfavorable angles between the cameras and pedestrian movements.

The simulation environment is a simplified urban setting—a flat, open area measuring $300 \times 160$ square feet. We assume there are no obstacles such as buildings or streets; pedestrians can move freely within this space.

The three PTZ cameras are positioned along the lower side of the environment at a height of 50 feet. Each camera has a 90-degree field of view and an aspect ratio of 1. The cameras have wide ranges for pan angles (from $-180^\circ$ to $180^\circ$), tilt angles (from $-90^\circ$ to $0^\circ$), and zoom ratios (up to 10x was utilized). For simplicity, we assume that the time required for a camera to adjust its settings (pan, tilt, and zoom) to a new configuration is constant, regardless of the previous and new parameters. This transition time is set to 1 second per camera. When a camera interrogates a pedestrian, it spends approximately 2 seconds capturing video, resulting in a total task duration of about 3 seconds per interrogation.

Pedestrians enter the environment from the north according to a Poisson process and may exit from any side after some time. This means the system has a limited window to capture video of each pedestrian before they leave the scene. The state of a pedestrian at time step $t$ is defined as $X(t) = [x(t), y(t), \dot{x}(t), \dot{y}(t)]^\top$, where $(x(t), y(t))$ represents the position and $(\dot{x}(t), \dot{y}(t))$ represents the velocity. At each frame, we obtain measurements of the pedestrian's features, forming an observation vector $Z(t) = [x_c(t), y_c(t), h(t), w(t)]^\top$, where $(x_c(t), y_c(t))$ is the center of the pedestrian's bounding box, and $(h(t), w(t))$ are its height and width.

We model the pedestrians' motion using a linear discrete-time dynamic system:
\begin{equation} \label{eq:objectDynamics}
    X(t+1) = AX(t) + W(t)
\end{equation}
and the observation model is given by:
\begin{equation} \label{eq:observationDynamics}
    Z(t+1) = Z(t) + V(t)
\end{equation}
In \eqref{eq:objectDynamics}, $A$ is the state transition matrix. $W(t)$ and $V(t)$ are zero-mean Gaussian noise vectors with small variances, representing process and measurement noise, respectively.\\
Cameras detect objects in both zoomed-in and zoomed-out modes, providing observations for each pedestrian. Using this information, we create a network flow model with predicted locations of pedestrians and solve the optimization problem to determine the cameras' schedules over the next $H$ periods. These schedules are stored in the cameras' memory, and the cameras execute the planned tasks period by period until a new plan is generated.
\subsection{Evaluation Metrics} \label{sec:metrics}
To evaluate and compare the performance of the surveillance systems, we use three key metrics:
\begin{itemize}
    \item Ratio of Covered Objects: This metric measures the proportion of pedestrians captured by the cameras, calculated as the number of interrogated pedestrians divided by the total number present in the environment.
    \item Average Wait Time: This metric represents the average time between a pedestrian's appearance in the environment and their interrogation by a camera, reflecting system responsiveness.
    \item Ratio of Lost Objects: This metric measures the proportion of pedestrians not captured before exiting the environment, calculated as the number of un-interrogated pedestrians divided by the total number.
\end{itemize}
These metrics provide a comprehensive evaluation of the systems' ability to capture pedestrians efficiently and promptly.
\subsection{Results}
\label{sec:resutls}
We conducted several simulation scenarios with varying parameters to evaluate the performance of our proposed systems. The key parameters we varied were the Poisson rate at which pedestrians enter the environment and the total number of pedestrians. These variations are critical for demonstrating the superiority of our method over alternatives like the master-slave configuration, especially under different crowd densities. In all scenarios, we compare three systems:
\begin{itemize}
    \item Flexible System: Our proposed system without group tracking.
    \item Flexible System with Group Tracking: Our proposed system enhanced with the group-tracking modification.
    \item Master-Slave System: A conventional system where one of the PTZ cameras is replaced with a wide-view static camera.
\end{itemize}
Below, we present the results for two representative scenarios.
\subsubsection{Scenario 1}
\begin{figure*}[t] % or [htbp] for more options
    \centering
    \begin{subfigure}{0.3\textwidth}
        \centering
        \includegraphics[width=\linewidth]{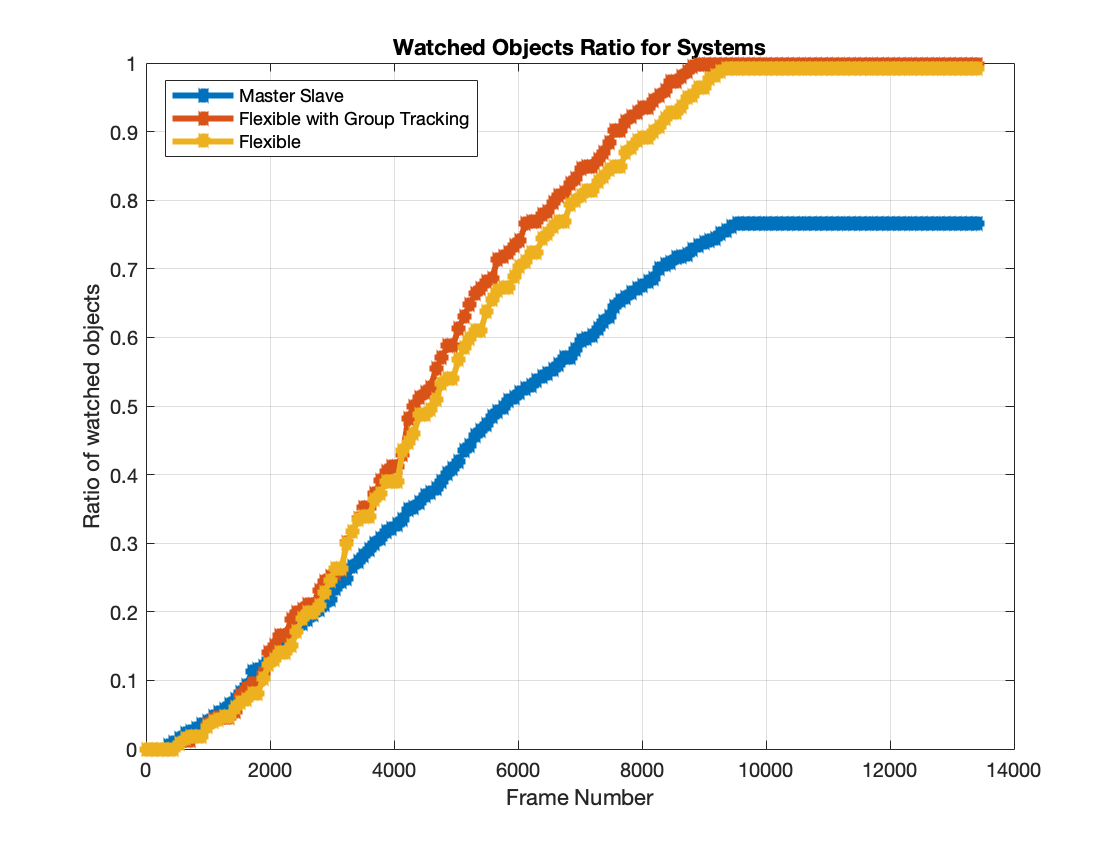}
        \caption{}
    \end{subfigure}
    \begin{subfigure}{0.3\textwidth}
        \centering
        \includegraphics[width=\linewidth]{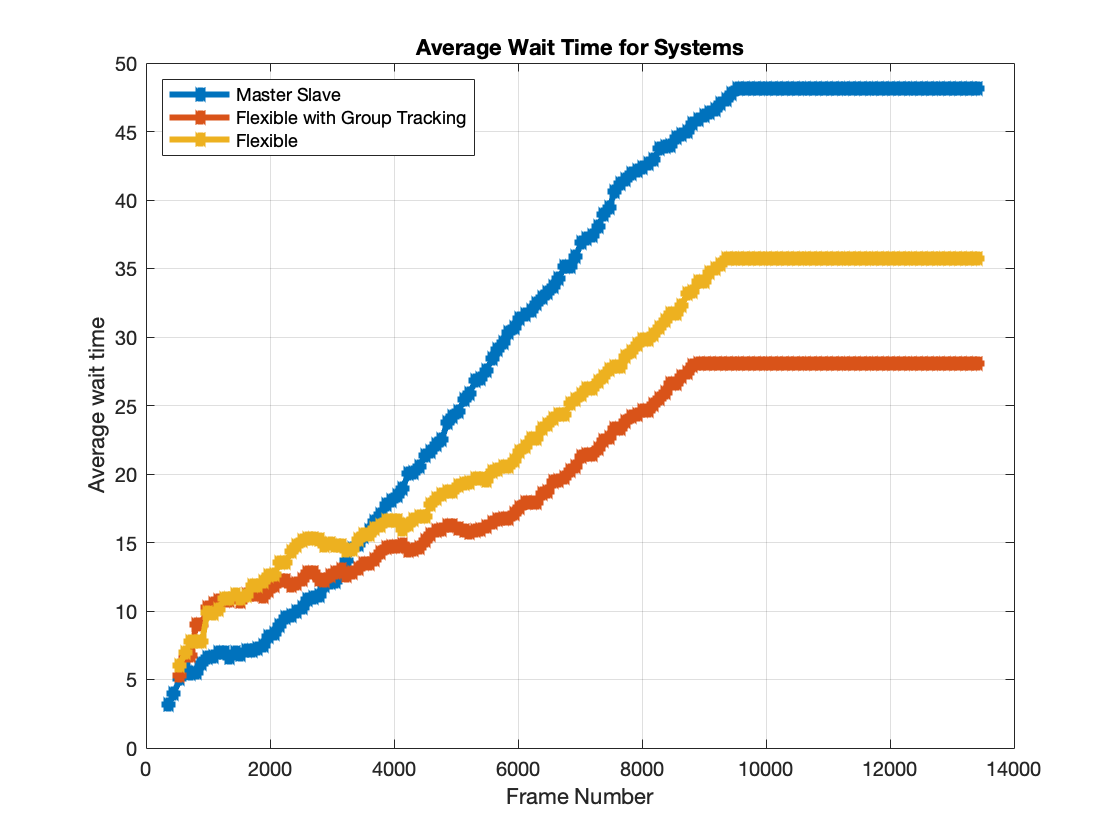}
        \caption{}
    \end{subfigure}
    \begin{subfigure}{0.3\textwidth}
        \centering
        \includegraphics[width=\linewidth]{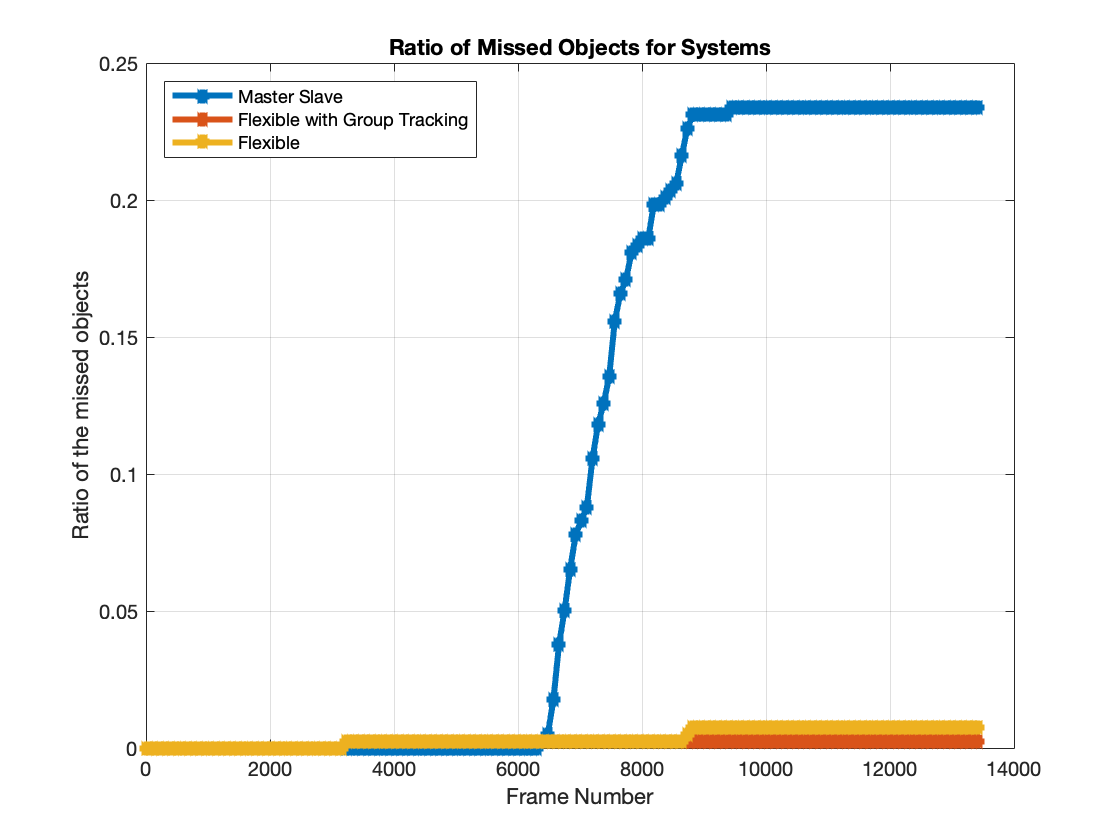}
        \caption{}
    \end{subfigure}
    % \vspace{4.5cm}
    \caption{Performance of different systems in a scenario with 400 pedestrians and pedestrian generation with Poisson rate of 1/20 per frame. (a) and (c) show the ratio of covered and missed objects which must add up to one. We can see that the flexible system with group tracking has a ratio of covered objects close to 1 (only one lost object for this system while the flexible system lost 3), and this metric is 0.75 for the master-slave system. (b) shows the average wait time for the systems.}
    \label{fig:p20n400}
\end{figure*}

In this scenario, pedestrians enter the environment at a Poisson rate of 1/20 per frame, and a total of 400 pedestrians appear over 450 seconds. Fig. \ref{fig:p20n400} shows the performance metrics for each system.

From the results, we observe that the flexible system with group tracking significantly outperforms both the standard flexible system and the master-slave system. Specifically:
\begin{itemize}
    \item Ratio of Covered Objects: The flexible system with group tracking achieves a coverage ratio near 1, indicating nearly all pedestrians were interrogated, while the master-slave system covers about 75 percent.
    \item Average Wait Time: The flexible system with group tracking demonstrates shorter wait times, indicating faster response to new pedestrians entering the environment.
    \item Ratio of Lost Objects: The flexible system with group tracking loses very few pedestrians (only one in this scenario), whereas the other systems miss more.
\end{itemize}

\subsubsection{Scenario 2}
\begin{figure*}[t] % or [htbp] for more options
    \centering
    \begin{subfigure}{0.3\textwidth}
        \centering
        \includegraphics[width=\linewidth]{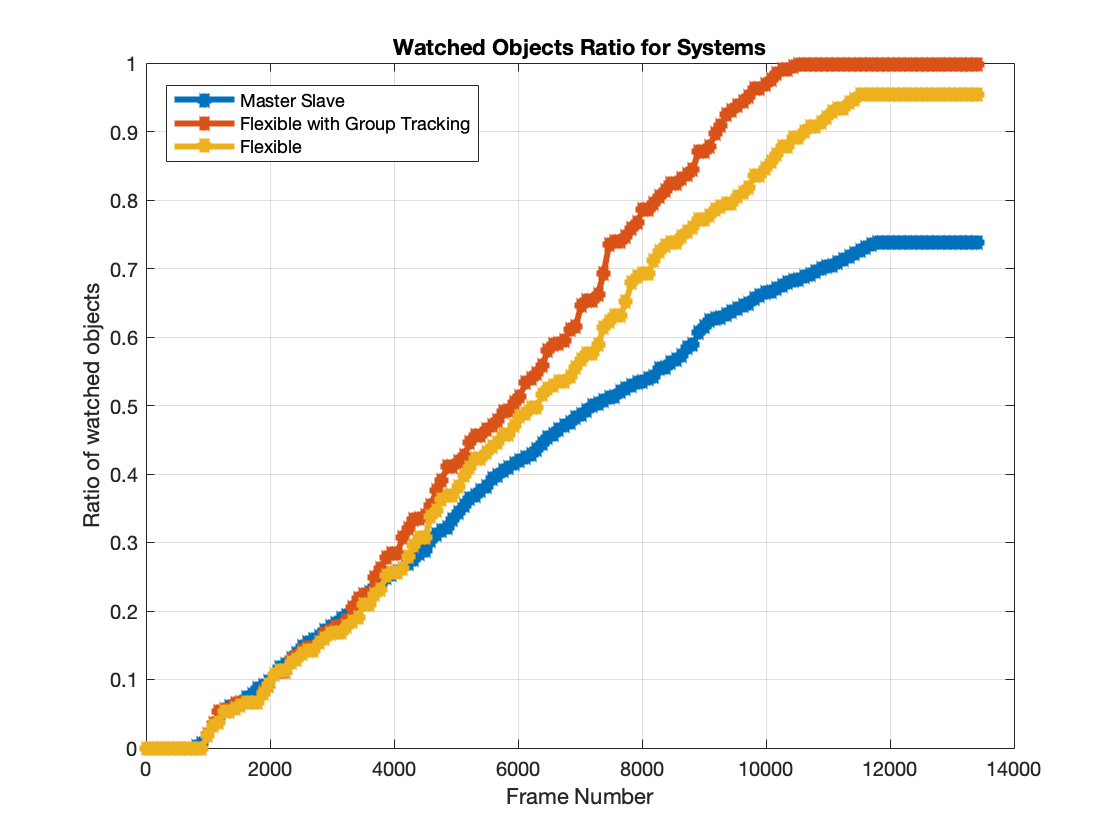}
        \caption{}
    \end{subfigure}
    \begin{subfigure}{0.3\textwidth}
        \centering
        \includegraphics[width=\linewidth]{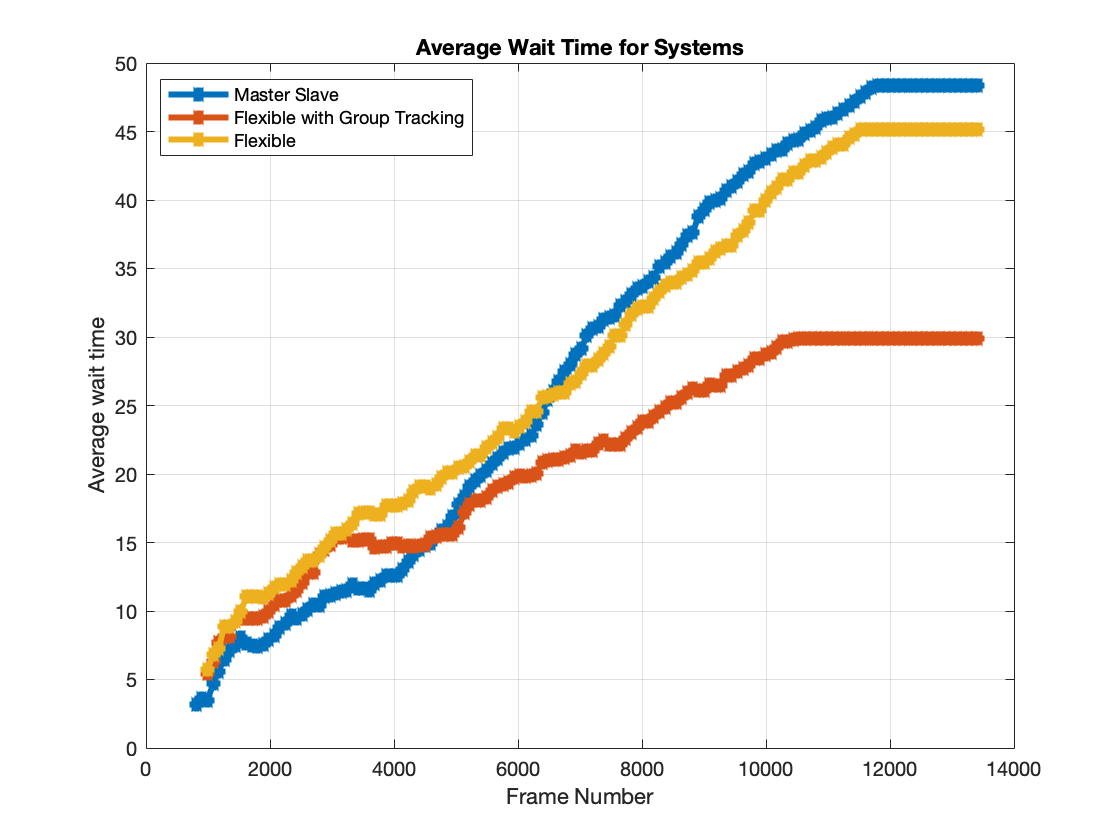}
        \caption{}
    \end{subfigure}
    \begin{subfigure}{0.3\textwidth}
        \centering
        \includegraphics[width=\linewidth]{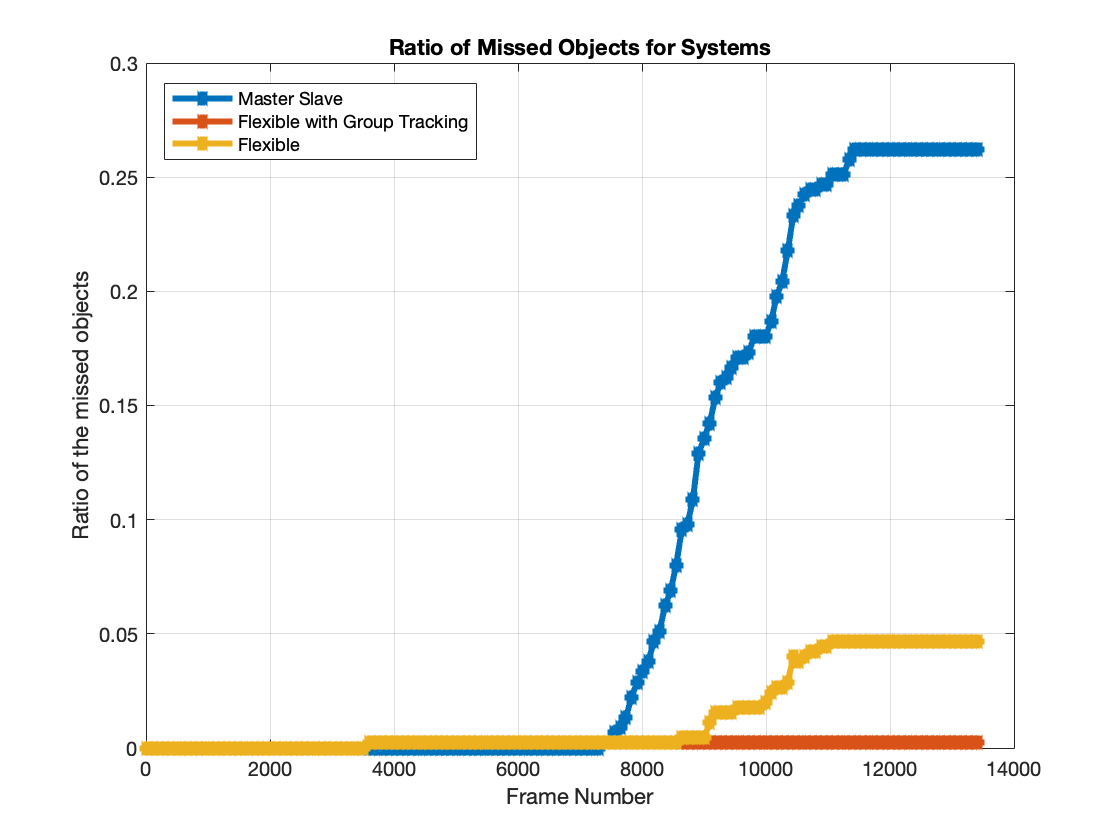}
        \caption{}
    \end{subfigure}
    % \vspace{4.5cm}
    \caption{Performance of different systems in a scenario with 450 pedestrians and pedestrian generation with Poisson rate of 1/18 per frame. Only one object is lost for the flexible system with group tracking while the performance of the other systems dropped dramatically.}
    \label{fig:p18n450}
\end{figure*}

In this more challenging scenario, pedestrians enter the environment at a Poisson rate of 1/18 per frame, and a total of 450 pedestrians appear over the simulation period. Fig. \ref{fig:p18n450} presents the performance metrics. We can observe that the flexible system with group tracking is still better than the others in all metrics.

Tables \ref{table1} and \ref{table2} present the quantitative results for two scenarios that we created and evaluated using our methods, demonstrating that both of our methods outperform the master-slave configuration.
\begin{table}[t]
\caption{Comparison of different methods across metrics for scenario 1}
\begin{center}
\begin{tabular}{|c|c|c|c|}
\hline
\textbf{Method}&\multicolumn{3}{|c|}{\textbf{Metrics}} \\
\cline{2-4} 
\textbf{Name} & \textbf{\textit{Watched Ratio}}& \textbf{\textit{\makecell{Average Wait \\ Time(s)}}}& \textbf{\textit{\makecell{Missed Objects \\ Ratio}}} \\
\hline
Master-Slave & 0.7663 & 48.10 & 0.2334 \\
\hline
Flexible & 0.9924 & 35.70 & 0.0075 \\
\hline
\makecell{Flexible with \\ group-tracking} & \textbf{0.9975} & \textbf{28.05} & \textbf{0.0025} \\
\hline
\end{tabular}
\label{table1}
\end{center}
\end{table}

\begin{table}[t]
\caption{Comparison of different methods across metrics for scenario 2}
\begin{center}
\begin{tabular}{|c|c|c|c|}
\hline
\textbf{Method}&\multicolumn{3}{|c|}{\textbf{Metrics}} \\
\cline{2-4} 
\textbf{Name} & \textbf{\textit{Watched Ratio}}& \textbf{\textit{\makecell{Average Wait \\ Time(s)}}}& \textbf{\textit{\makecell{Missed Objects \\ Ratio}}} \\
\hline
Master-Slave & 0.7378 & 48.35 & 0.2622 \\
\hline
Flexible & 0.9533 & 45.18 & 0.0467 \\
\hline
\makecell{Flexible with \\ group-tracking} & \textbf{0.9978} & \textbf{29.9} & \textbf{0.0022} \\
\hline
\end{tabular}
\label{table2}
\end{center}
\end{table}
\section{Discussion}
The results from these scenarios illustrate the effectiveness of our proposed method, particularly when enhanced with group tracking. The group-tracking modification significantly improves performance by allowing the system to capture multiple pedestrians simultaneously, thus optimizing camera utilization. Compared to the master-slave system, our flexible system with group tracking demonstrates higher coverage, capturing videos of nearly all pedestrians even as the crowd density increases. It also exhibits lower wait times, responding more quickly to new pedestrians and reducing the delay between their arrival and interrogation. Furthermore, the system maintains scalability, showing high performance levels even under more challenging conditions with more pedestrians and higher arrival rates.

These improvements validate the advantages of our network flow optimization approach combined with group tracking, confirming its potential for real-world surveillance applications where efficiency and scalability are critical. As an example, for the initial scenario, utilizing MATLAB's linprog on an Apple M2 Pro CPU with 32GB of RAM, the peak computation time for creating a new camera planning was 0.0509 seconds, recorded for the largest evaluated graph comprising 367 vertices and 1326 arcs.

\section{Conclusion}
In this paper, we presented a novel scheduling algorithm for smart video capture using PTZ cameras in surveillance systems. By utilizing Kalman filters to track objects within the environment, we accurately predicted future object locations and constructed a network flow model to optimize the scheduling of camera tasks. Our simulations demonstrated that this method outperforms the traditional master-slave camera network configuration across key metrics such as the ratio of covered objects and average wait time.

An important enhancement to our system was the use of a greedy set cover algorithm to introduce group-tracking nodes. By grouping nearby targets based on predicted locations, we improved performance, reduced computational complexity, minimized redundant tracking, and enabled simultaneous multi-object capture.

While our experimental results showcase the effectiveness of our algorithm in terms of the number of captured objects and reduced waiting times, there are opportunities for further improvement. One direction for future work is to relax the assumption that all cameras have the same transition time between tasks. In real-world scenarios, camera adjustment times for pan, tilt, and zoom depend on the initial and final settings. Incorporating variable transition times based on the specific PTZ adjustments required for each task would make the model more accurate and could lead to more efficient scheduling. Additionally, implementing the system in practical settings will allow us to address challenges such as processing speed, communication delays, and varying environmental conditions. By extending our research in these directions and incorporating variable transition times, we aim to develop a robust and scalable solution for intelligent surveillance that can adapt to dynamic and complex environments.

\bibliographystyle{IEEEtran}
\bibliography{sources.bib}

\end{document}